\documentclass[a4paper,reqno,twoside]{amsart}

\usepackage[lite,initials]{amsrefs}
\usepackage{amsthm,amsmath,amssymb}

\usepackage{paralist}

\newtheorem{theorem}{Theorem}
\newtheorem{proposition}[theorem]{Proposition}

\theoremstyle{definition}

\numberwithin{theorem}{section}
\numberwithin{equation}{section}

\providecommand{\norm}[1]{\left\lVert#1\right\rVert}
\providecommand{\Id}{\operatorname{Id}}

\begin{document}

\title{Subspaces of almost Daugavet spaces}
\author{Simon L\"ucking}
\address{Departement of Mathematics, Freie Universit\"at Berlin, Arnimallee 6,
					14195 Berlin, Germany}
\email{simon.luecking@fu-berlin.de}

\begin{abstract}
	We study the almost Daugavet property, a generalization of the Daugavet
	property. It is analysed what kind of subspaces and sums of Banach
	spaces with the almost Daugavet property have this property as well. The
	main result of the paper is: if $Z$ is a closed subspace of a separable
	almost Daugavet space $X$ such that the quotient space $X/Z$ contains no
	copy of	$\ell_1$, then $Z$ has the almost Daugavet property, too.	
\end{abstract}

\maketitle

\section{Introduction}

V. Kadets, V. Shepelska, and D. Werner characterized in
\cite{Wernerthickness} a class of separable Banach spaces in three very
different ways.

\begin{theorem}
	\label{theoremaequivalenz}
	For a separable Banach space $X$ the following conditions are equivalent:
	\begin{enumerate}[\upshape(i)]
		\item $S_X$ does not admit a finite $\varepsilon$-net consisting of
			elements of $S_X$ for any $\varepsilon < 2$.
		\item there is a sequence $\left(e_n\right)_{n \in \mathbb{N}}$ in $B_X$
			such that for every $x \in X$
			\begin{equation*}
				\lim_{n \rightarrow \infty} \norm{x+e_n} = \norm{x} +1.
			\end{equation*}
		\item there is a norming subspace $Y \subset X^*$ such that the equation
			\begin{equation}
				\label{Daugavetequation}
				\norm{\Id + T} = 1 + \norm{T}
			\end{equation}
			holds true for every rank-one operator $T: X \rightarrow X$ of the form
			$T= y^*\otimes x$, where $x\in X$ and $y^* \in Y$.
	\end{enumerate}	
\end{theorem}

Recall that a subspace $Y \subset X^*$ is said to be norming if for every
$x \in X$
\begin{equation*}
	\sup_{y^* \in S_Y} \left|y^*(x)\right| = \norm{x}.
\end{equation*}

Condition (i) can be rephrased using the so-called \emph{thickness} of a 
Banach space $X$, that was introduced by R. Whitley in \cite{Whitley}.
We call a set $F$ an inner $\varepsilon$-net of $S_X$ if $F \subset S_X$ 
and for every $x \in S_X$, there is an $y \in F$ with $\norm{x-y} \leq
\varepsilon$. Then the thickness $T(X)$ of a Banach space $X$ is defined by
\begin{equation*}
	T(X) = \inf \left\{ \varepsilon>0 \colon \text{there exists a finite
	inner $\varepsilon$-net of $S_X$} \right\}.
\end{equation*}
Essentially, the thickness is a kind of inner measure of non-compactness
of the unit sphere $S_X$.\\
R. Whitley showed that $1 \leq T(X) \leq 2$ if $X$ is infinite-dimensional, 
in particular that $T(l_p)=2^{1/p}$ for $1 \leq p < \infty$ and $T(C(K))=2$
if $K$ has no isolated points.\\
So condition (i) means that $T(X) = 2$.

Condition (ii) links the investigation to the theory of types, that were used
by \mbox{J.-L.} Krivine and B. Maurey in \cite{KrivineMaurey}. A type on a
separable Banach space is a function of the form
\begin{equation*}
	\tau(x) = \lim_{n \rightarrow \infty} \norm{x + e_n}
\end{equation*}
for some bounded sequence $\left(e_n\right)_{n \in \mathbb{N}}$. In
\cite{KrivineMaurey} a convolution of types is defined and by means of this 
convolution the notion of $\ell_1$-types is introduced. A special $\ell_1$-type
is generated by a sequence $\left(e_n\right)_{n \in \mathbb{N}}$ satisfying
\begin{equation*}
	\tau(x)= \lim_{n \rightarrow \infty} \norm{x + e_n} = \norm{x} +1.
\end{equation*}
We call a sequence $\left(e_n\right)_{n \in \mathbb{N}}$ in $B_X$ satisfying
this equation a \emph{canonical $\ell_1$-type sequence}. The following
proposition shows us that a canonical $\ell_1$-type sequence (or a suitable
subsequence) almost behaves like the canonical base of $\ell_1$.

\begin{proposition}
	\label{teilfolgekanonischerl1typ}
	Let $X$ be a separable Banach space, $E_1 \subset E_2 \subset E_3 \subset
	\ldots$ a sequence of finite-dimensional subspaces whose union is dense in
	$X$, and $\left(\varepsilon_n\right)_{n \in\mathbb{N}}$ a	sequence of
	positive real numbers that converges monotonically towards zero. 
	If $\left( e_n \right)_{n \in \mathbb{N}}$ is a canonical $\ell_1$-type
	sequence then we can pass to a subsequence $\left( f_n \right)_
	{n \in \mathbb{N}}$ of $\left( e_n \right)_{n \in \mathbb{N}}$ such that for
	every $x \in E_n$ and	$y = \sum_{k=n+1}^{M}\lambda_k f_k$
	\begin{equation*}
		\norm{x+y} \geq \left( 1 - \varepsilon_n \right) \norm{x}
		+ \sum_{k=n+1}^{M} \left(1-\varepsilon_{k-1} \right) 
		\left|\lambda_k \right|.
	\end{equation*}
	\begin{proof}
		Fix a sequence of positive real numbers $\left( \delta_n \right)_{n \in
		\mathbb{N}}$ with $\prod_{k=n}^\infty \left( 1- \delta_k \right) \geq 1-
		\varepsilon_n$ for all $n \in \mathbb{N}$. There is a subsequence 
		$\left( f_n \right)_{n \in \mathbb{N}}$ of $\left( e_n \right)_{n \in
		\mathbb{N}}$ such that for every $x \in \operatorname{span}(E_n \cup
		\{f_1, \ldots, f_n\})$ and $\lambda \in \mathbb{K}$ 
		\begin{equation*}
			\norm{x+\lambda f_{n+1}} \geq \left(1-\delta_n\right)
			\left(\norm{x}+ \left|\lambda\right|\right). 
		\end{equation*}
		Then we have for every $x \in E_n$ and $y = \sum_{k=n+1}^{M}{\lambda_k
		f_k}$ as a consequence of the choice of $\left( \delta_n \right)_{n \in
		\mathbb{N}}$
		\begin{align*}
			\norm{x+y} = \norm{x+\sum_{k=n+1}^{M}\lambda_k f_k}
			&\geq \left(1-\delta_{M-1}\right) \norm{x+\sum_{k=n+1}^{M-1}
			 {\lambda_k f_k}}
			 +\left(1-\delta_{M-1}\right) \left|\lambda_M\right| \\
			&\geq \dotsb  
			 \geq \prod_{k=n}^{M-1}{(1-\delta_k)} \norm{x} + \sum_{k=n+1}^{M} 
			 \left(1-\delta_{k-1} \right) \left|\lambda_k \right| \\
			&\geq	
			 \left( 1 - \varepsilon_n \right) \norm{x}
			 + \sum_{k=n+1}^{M} \left(1-\varepsilon_{k-1} \right) 
			 \left|\lambda_k \right|. 
			 \qquad\qquad \text{ }\text{ }\mbox{\qedhere}
		\end{align*}
	\end{proof}	
\end{proposition}

Condition (iii) links the investigation to the theory of Banach spaces with the
\emph{Daugavet property} introduced by V. Kadets, R. Shvidkoy, G. Sirotkin,
and D. Werner in \cite{Daugavetproperty}. 
We say that a Banach space $X$ has the \emph{Daugavet property with respect to
Y} if the \emph{Daugavet equation} (\ref{Daugavetequation}) holds true for
every rank-one operator $T: X \rightarrow X$ of the form $T = y^* \otimes x$,
where $x \in X$ and $y^* \in Y$. A Banach space $X$ is called an \emph{almost
Daugavet space} or a space with the \emph{almost Daugavet property} if it has
the Daugavet property with respect to some norming subspace $Y \subset X^*$.
This definition, that was introduced in \cite{Quotients}, is a generalization
of the Daugavet property. Recall that a Banach space $X$ has the Daugavet
property if it has the Daugavet property with respect to $X^*$.

In the following section we study some structural properties of the class of
almost Daugavet spaces, and our main tool will be theorem
\ref{theoremaequivalenz}.
We show that the $\ell_1$- and $\ell_\infty$-sum of separable almost Daugavet
spaces are almost Daugavet spaces.\\
Concerning subspaces, we show results about $L$-summands and $M$-ideals. But 
the main result will be that the almost Daugavet property of a separable 
Banach space $X$ is inherited by a closed subspace $Z$ if the quotient space
$X/Z$ contains no copy of $\ell_1$.\\
Special thanks go to V. Kadets for giving an important hint to achieve this 
result in such a general setting.

\section{Subspaces and sums of almost Daugavet spaces}

\subsection{Sums}

The almost Daugavet property depends very much on the norm of a space.
So it can't be expected that arbitrary sums of almost Daugavet spaces
still have the almost Daugavet property. But there are the following results
concerning $\ell_1$- and $\ell_\infty$-sums.

\begin{proposition}
	If $X$ and $Y$ are two separable almost Daugavet spaces, then 
	$X \oplus_\infty Y$ is an almost Daugavet space as well.
	\begin{proof}
		As a consequence of theorem \ref{theoremaequivalenz}, $X$ as well as $Y$
		contains a canonical $\ell_1$-type sequence $\left(e_n\right)_{n\in
		\mathbb{N}}$ resp. $\left(f_n\right)_{n\in \mathbb{N}}$. Every member
		of $\left(\left(e_n,f_n\right)\right)_{n \in \mathbb{N}}$ is an 
		element of $B_{X\oplus_\infty Y}$ and for $(x,y) \in X \oplus_\infty Y$
		we have
		\begin{align*}
			\lim_{n \rightarrow \infty} \norm{(x,y)+(e_n,f_n)} 
			&=\lim_{n \rightarrow \infty}
			  \max\left\{\norm{x+e_n},\norm{y+f_n}\right\}\\
			&=\max\left\{ \lim_{n \rightarrow \infty}\norm{x+e_n},
				\lim_{n \rightarrow \infty}\norm{y+f_n}\right\}\\
			&=\max\left\{\norm{x}+1,\norm{y}+1\right\}\\
			&=\max\left\{\norm{x},\norm{y}\right\}+1=\norm{(x,y)}+1.			 
		\end{align*}
		So $\left(\left(e_n,f_n\right)\right)_{n \in \mathbb{N}}$ is a 
		canonical $\ell_1$-type sequence in $X \oplus_\infty Y$ and 
		$X \oplus_\infty Y$ is an almost Daugavet space by
		theorem \ref{theoremaequivalenz}.
	\end{proof}
\end{proposition}

\begin{proposition}
	Let $X$ and $Y$ be separable Banach spaces. If $X$ has the almost
	Daugavet property, then so does $X \oplus_1 Y$.
	\begin{proof}
		Because of theorem \ref{theoremaequivalenz}, $X$ contains a canonical
		$\ell_1$-type sequence $\left(e_n\right)_{n\in \mathbb{N}}$.
		For $(x,y) \in X \oplus_1 Y$ we get
		\begin{align*}
			\lim_{n \rightarrow \infty} \norm{(x,y)+(e_n,0)} 
			&= \lim_{n \rightarrow \infty} \left(\norm{x+e_n} + \norm{y}\right)
			= \lim_{n \rightarrow \infty} \norm{x+e_n} + \norm{y}\\
			&= \norm{x} + \norm{y} + 1 = \norm{(x,y)} + 1,
		\end{align*}
		and thus $\left(\left(e_n,0\right)\right)_{n \in \mathbb{N}}$ is a
		canonical $\ell_1$-type sequence. Therefore $X \oplus_1 Y$ is an almost 
		Daugavet space as a result of theorem \ref{theoremaequivalenz}.
	\end{proof}
\end{proposition}

\subsection{Subspaces}

Now we are interested in the question which subspaces of almost Daugavet spaces
inherit the almost Daugavet property. Since the almost Daugavet property 
depends very much of the norm of a space, we are going to study subspaces
fulfilling adequat conditions.

\begin{proposition}
	Let $X$ be a separable almost Daugavet space, $J$ an $L$-summand and
	$\widehat{J}$ the
	complementary $L$-summand (i.e., $X=J \oplus_1 \widehat{J}$). Then $J$ or
	$\widehat{J}$ has the almost Daugavet property.
	\begin{proof}
		As a result of theorem \ref{theoremaequivalenz}, it suffices to show that
		$T(J) =2$ or $T(\widehat{J})=2$.\\
		Let us suppose that $T(J) < 2$ and $T(\widehat{J}) < 2$. Then there exist
		for an $\varepsilon > 0$ inner $(2-2\varepsilon)$-nets of $S_J$ and
		$S_{\widehat{J}}$. Let $\left\{ j_1, \ldots, j_l\right\}$ be such an inner
		$(2-2\varepsilon)$-net of $S_J$ and $\{ \widehat{j}_1, \ldots,
		\widehat{j}_m\}$ one of $S_{\widehat{J}}$. Furthermore, let $\left\{
		\lambda_1, \ldots, \lambda_n\right\}$ be an inner $(\varepsilon/2)$-net of
		$[0,1]$.\\
		As $J$ and $\widehat{J}$ are $L$-summands, the unit sphere of $X$ can be
		expressed by
		\begin{equation*}
			S_X = \left\{\left(\lambda j , (1-\lambda)\widehat{j}\right) \colon
			\lambda \in [0,1], j \in S_J, \widehat{j} \in S_{\widehat{J}} \right\}.
		\end{equation*}
		Given $x = (\lambda j , (1-\lambda)\widehat{j})	\in S_X$, we can pick 
		$\lambda_0 \in \left\{ \lambda_1, \ldots, \lambda_n\right\}$ with
		$\left|\lambda_0 - \lambda\right| \leq \varepsilon/2$, $j_0 \in \left\{
		j_1, \ldots, j_l\right\}$ with $\norm{j_0 - j} \leq	2-2 \varepsilon$ and
		$\widehat{j}_0 \in \{ \widehat{j}_1, \ldots,
		\widehat{j}_m\}$ with ${\|\widehat{j}_0 - \widehat{j}\| \leq
		2-2 \varepsilon}$. Then we get with $y = (\lambda_0 j_0 ,
		(1-\lambda_0)\widehat{j}_0)$
		\begin{align*}
			\norm{x- y} &=
			\norm{\lambda j - \lambda_0 j_0}
			+ \norm{(1-\lambda)\widehat{j} - (1- \lambda_0)\widehat{j}_0} \\
			&=\norm{(\lambda - \lambda_0)j +\lambda_0 (j-j_0)}
			+ \norm{(\lambda_0 - \lambda)\widehat{j}
			 + (1-\lambda_0)(\widehat{j}-\widehat{j}_0)}\\
			&\leq |\lambda - \lambda_0|\left(\norm{j}+\|\widehat{j}\|\right)
				+\left(\lambda_0 \norm{j-j_0} + (1-\lambda_0)\|\widehat{j} -
				\widehat{j}_0\|\right)\\
			&\leq \varepsilon + \left(2 - 2 \varepsilon\right) = 2 - \varepsilon.	
		\end{align*}
		Thus $\left\{(\lambda j , (1-\lambda)\widehat{j}) \colon
		\lambda \in \left\{\lambda_1, \ldots, \lambda_n\right\}, j \in 
		\left\{j_1, \ldots, j_l\right\}, \widehat{j} \in
		\{\widehat{j}_1, \ldots, \widehat{j}_m\}\right\}$ is an inner
		$(2-\varepsilon)$-net of $S_X$. But this is because of theorem 
		\ref{theoremaequivalenz} a contradiction to the	assumption that $X$ has
		the almost Daugavet property. 
	\end{proof}
\end{proposition}

Recall that an $M$-ideal in a Banach space $X$ is a closed subspace $J$ such
that $X^*$ decomposes as $X^* = Y \oplus_1 J^\bot$ for some closed subspace
$Y$ of $X^*$, where ${J^\bot = \left\{x^* \in X^* \colon x^*(x) = 0 
\text{ for all $x \in J$}\right\}}$.\\
In \cite{Daugavetproperty}*{Proposition 2.10.} it is shown that $M$-ideals of
Banach spaces with the Daugavet property have the Daugavet property as well.
The same proof works for almost Daugavet spaces.

\begin{proposition}
	Let $X$ be a Banach space, $Y$ a norming subspace of $X^*$ and suppose that 
	$X$ has the Daugavet property with respect to $Y$. If $J \neq \{0\}$ is an
	$M$-ideal in $X$, then $J$ has the Daugavet property with respect to
	$\left\{y^*|_J \colon	y^* \in Y \right\}$. So $J$ is an almost Daugavet
	space as well.
\end{proposition}

On the other hand, almost Daugavet spaces are in a certain sense ``big''. They
have thickness 2 and contain a copy of $\ell_1$. So in the next theorem we
consider ``big'' subspaces of an almost Daugavet space.

\begin{theorem}
	Let $X$ be a separable almost Daugavet space and $Z$ a closed subspace of
	$X$. If the quotient space $X/Z$ contains no copy of $\ell_1$, then $Z$ has
	the almost Daugavet property, too.
	\begin{proof}
		By theorem \ref{theoremaequivalenz}, it suffices to construct a
		canonical $\ell_1$-type sequence in $Z$.\\ 
		Let $E_1 \subset E_2 \subset E_3 \subset \ldots$ be a sequence of 
		finite-dimensional subspaces whose union is dense in $X$. Also, fix a
		sequence of positive real numbers $\left(\varepsilon_n\right)
		_{n \in\mathbb{N}}$ that converges monotonically towards zero. As a
		consequence of theorem \ref{theoremaequivalenz} and proposition
		\ref{teilfolgekanonischerl1typ}, there is a canonical $\ell_1$-type
		sequence $\left(e_n\right)
		_{n \in \mathbb{N}}$ in $B_X$ satisfying the following condition: for
		every $x \in E_n$ and every $y = \sum_{k=n+1}^M{\lambda_k e_k}$ we have
		\begin{equation}
		  \label{KodimensionI}
			\norm{x+y} \geq \left( 1 - \varepsilon_n \right) \norm{x}
		  + \sum_{k=n+1}^{M}\left(1-\varepsilon_{k-1} \right)\left|\lambda_k
		  \right|.		
		\end{equation}
		So the space $\overline{\operatorname{span}}\left\{e_n \colon n \geq k
		\right\}$ is
		isomorphic to $\ell_1$ for each $k \in \mathbb{N}$. Since $X/Z$ contains
		no copy of $\ell_1$, the quotient map $\pi: X \rightarrow X/Z$ fails to 
		be bounded below on each subspace $\overline{\operatorname{span}}\left\{e_n
		\colon n \geq k \right\}$, i.e., for each $k \in \mathbb{N}$, there is no
		$C > 0$ with $C \norm{x} \leq \norm{\pi(x)}$ for all $x \in
		\overline{\operatorname{span}}\left\{e_n	\colon n \geq k\right\}$.\\
		As a result of this, there is a linear combination
		$\sum_{k=1}^{m_1}\lambda_ke_k$ with $\norm{\sum_{k=1}^{m_1}
		\lambda_ke_k}=1$ and $\norm{\pi(\sum_{k=1}^{m_1}
		\lambda_ke_k)} \leq \varepsilon_1/4$. 
		By the definition of the quotient norm, there is an $\hat{f_1} \in Z$ with
		$\|\hat{f_1}-\sum_{k=1}^{m_1}\lambda_ke_k\| \leq \varepsilon_1/2$, and we
		get	$\norm{f_1-\sum_{k=1}^{m_1}\lambda_ke_k} \leq \varepsilon_1$ 
		with $f_1 = \hat{f_1}/\|\hat{f_1}\|$. Going on like this, we get for each
		$n \in \mathbb{N}$ a linear combination
		$\sum_{k=m_{n-1}+1}^{m_n}\lambda_ke_k$ with
		$\|\sum_{k=m_{n-1}+1}^{m_n}\lambda_ke_k\|=1$ and an element $f_n \in S_Z$
		with ${\|f_n - \sum_{k=m_{n-1}+1}^{m_n}\lambda_ke_k\| \leq
		\varepsilon_n}$.\\
		It remains to show that $\left(f_n\right)_{n \in \mathbb{N}}$ is a
		canonical $\ell_1$-type sequence. Given	$x \in X$ and $\varepsilon > 0$,
		we can assume without loss of generality that ${x \in E_{\hat{n}}}$ for an
		$\hat{n} \in \mathbb{N}$,	because the union of all subspaces $E_n$ is
		dense. The sequence $\left(\varepsilon_n\right)_{n \in\mathbb{N}}$
		converges monotonically towards zero. So there is an
		$N_0 \in \mathbb{N}$ with ${\varepsilon_{n-1} \norm{x} + \varepsilon_{n-1}
		+\varepsilon_n \leq \varepsilon}$ for all $n > N_0$. Then we get for all 
		$n > \max\left\{\hat{n},N_0\right\}$ using (\ref{KodimensionI}) and 
		$m_{n-1}\geq n-1$
		\begin{align*}
			\norm{x} + 1 \geq \norm{x + f_n} &\geq
			\norm{x + \sum_{k=m_{n-1}+1}^{m_n}\lambda_ke_k} 
			- \norm{\sum_{k=m_{n-1}+1}^{m_n}\lambda_ke_k
			+f_n}\\
			&\geq (1-\varepsilon_{n-1})\norm{x} + \sum_{k=m_{n-1}+1}^{m_n}
			(1-\varepsilon_{k-1})|\lambda_k| - \varepsilon_n \\
			&\geq (1-\varepsilon_{n-1})\norm{x} + (1-\varepsilon_{n-1})
			\sum_{k=m_{n-1}+1}^{m_n}|\lambda_k| - \varepsilon_n\\
			&\geq (1-\varepsilon_{n-1})\norm{x} + 
			(1-\varepsilon_{n-1}) - \varepsilon_n\\
			&\geq \norm{x} +1 - \varepsilon.
		\end{align*}
		So $\lim_{n \rightarrow \infty} \norm{x + f_n} = \norm{x} + 1$ and
		$\left(f_n\right)_{n \in \mathbb{N}}$ is a canonical $\ell_1$-type
		sequence in $Z$.
		\end{proof}
\end{theorem}

\end{document}